\documentclass[11pt]{article}
\usepackage{epsfig}
\usepackage{amssymb}

 \newtheorem{theorem}{Theorem}
 
 \newtheorem{corollary}[theorem]{Corollary}

\newtheorem{conjecture}[theorem]{Conjecture}

\newenvironment{proof}{{\it Proof:\/}}{$\Box$\vskip 0.08in}

\title{Positive knots have negative signature\footnote{The paper was published in 
{\it Bull. Polish Acad. Sci.: Math.} 37, no. 7-12, 1989, 559-562.}}

\author{Jozef H. Przytycki}
\date{}
\begin{document}
\maketitle 
\date{}

{\bf Abstract.} We show that if a nontrivial link in $R^3$ has a diagram with all
crossings positive 
(\parbox{.9cm}{\psfig{figure=L+maly.eps,height=0.7cm}}) 
then the signature of the link is negative. It
settles the old folklore conjecture.\\
 \ \\ \ \\


It was asked by Birman, Williams, and Rudolph whether nontrivial
Lorentz knots \cite{B-W} have always positive signature. Lorentz
knots are examples of positive braids (in our convention they have
all crossings negative so they are negative links). It was
shown by Rudolph \cite{R} that positive braids have positive
signature (if they represent nontrivial links). Murasugi has shown
that nontrivial, alternating, positive links have negative
signature. \\
Here we solve the conjecture in general.

\begin{theorem}\label{thm1}
  Let $L$ be a nontrivial link which has a diagram with all
  crossings positive (i.e. L is positive), then the signature of
  $L$, $\sigma(L)<0$.
\end{theorem}
\begin{proof}
  Our main tool is the result of Murasugi (compare \cite{P}) which
  says that if two links $L_+$ and $L_-$ have identical diagrams
  except near one crossing where they look as on Fig. 1 
\centerline{\psfig{figure=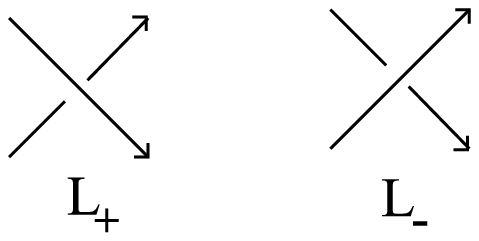}}
\centerline{Fig. 1.}
 then $\sigma(L_+) \leq \sigma(L_-)$. First
  assume that $L$ is a knot. Consider a positive diagram of
  $L$ (also denoted by $L$) with minimal number of crossings.
Consider an innermost $1$--gon in $L$. Now move along $L$ starting
from b (Fig. 2). Let $p_1$ be the last enter of $L$
into the $1$--gon. Now change the overcrossing to undercrossing 
in $L$ in such a way that, starting from $b$,
\ \\ \ \\
\centerline{\psfig{figure=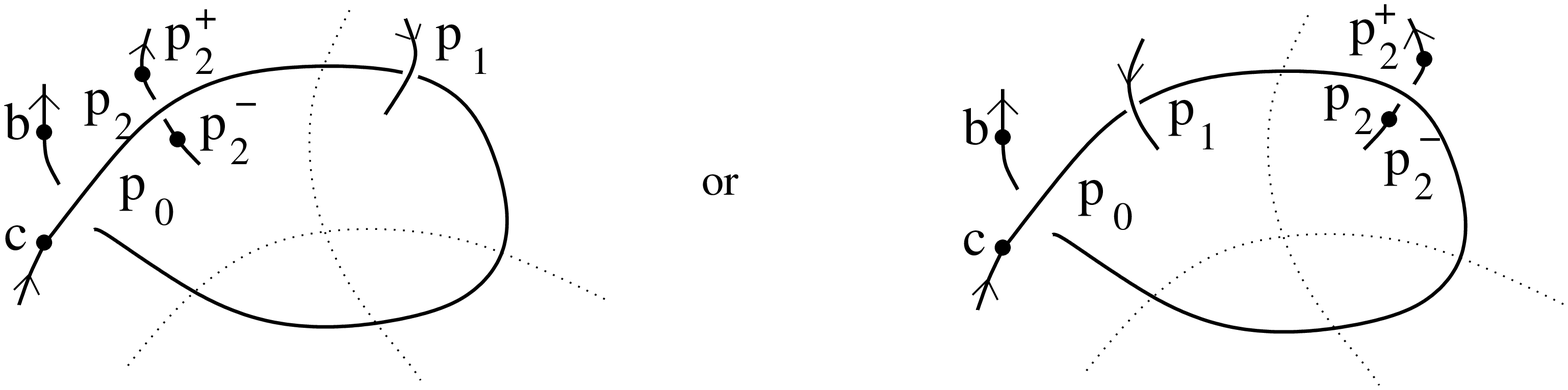,height=3.4cm}}
\centerline{Fig. 2.}
\ \\
the new diagram, $L'$, is descending except the crossing $p_2$. Therefore the
crossings $p_0$, $p_1$ and $p_2$ are not changed (Fig. 2). 
By the Murasugi result $\sigma(L') \geq \sigma(L)$.
We claim the $L'$ is a diagram of the right handed trefoil knot
(\parbox{1cm}{\psfig{figure=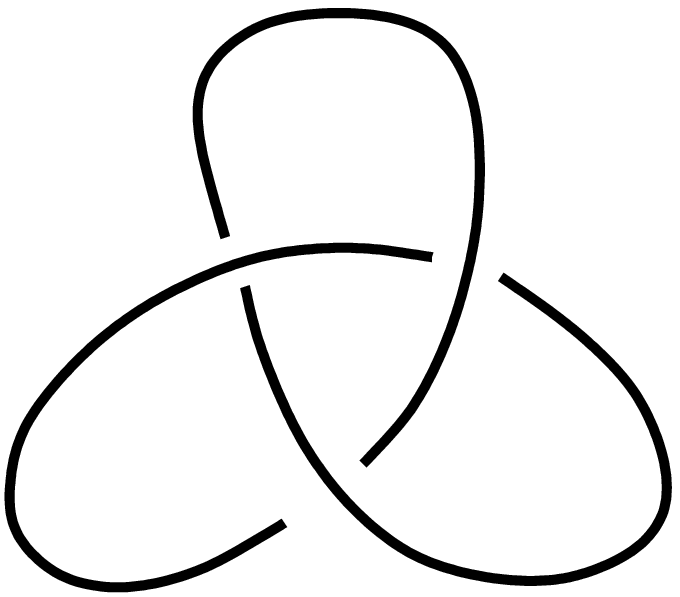, width=1cm}}) and
therefore $\sigma(L')=-2$, so $\sigma(L) < 0$. It remains to prove
the
above claim.\\
Let us assume that the point $b$ is on the level 3 and that the
diagram descends to level 2 just before $p_2$ (say in $p_2^-$)
then the underscrossing reach the  lever $-1$ and ascends back to
level 1 at c, finally reaching the level 0 just before b. From
this point of we ascend quickly to the level 3 at b. We can assume
that the 1-gon is convex. Now we can easily deform $L'$ by isotopy
so that the part $b$,$p_2^-$ is the straight line and the part
$p_2^+,c$, is a simple arc which is descending and disjoint from
other parts of the diagram (except ends). Therefore $L'$ is
isotopic to the diagram of Fig. 3 which represents the
right handed trefoil knot. \\
\ \\ \ \\
\centerline{\psfig{figure=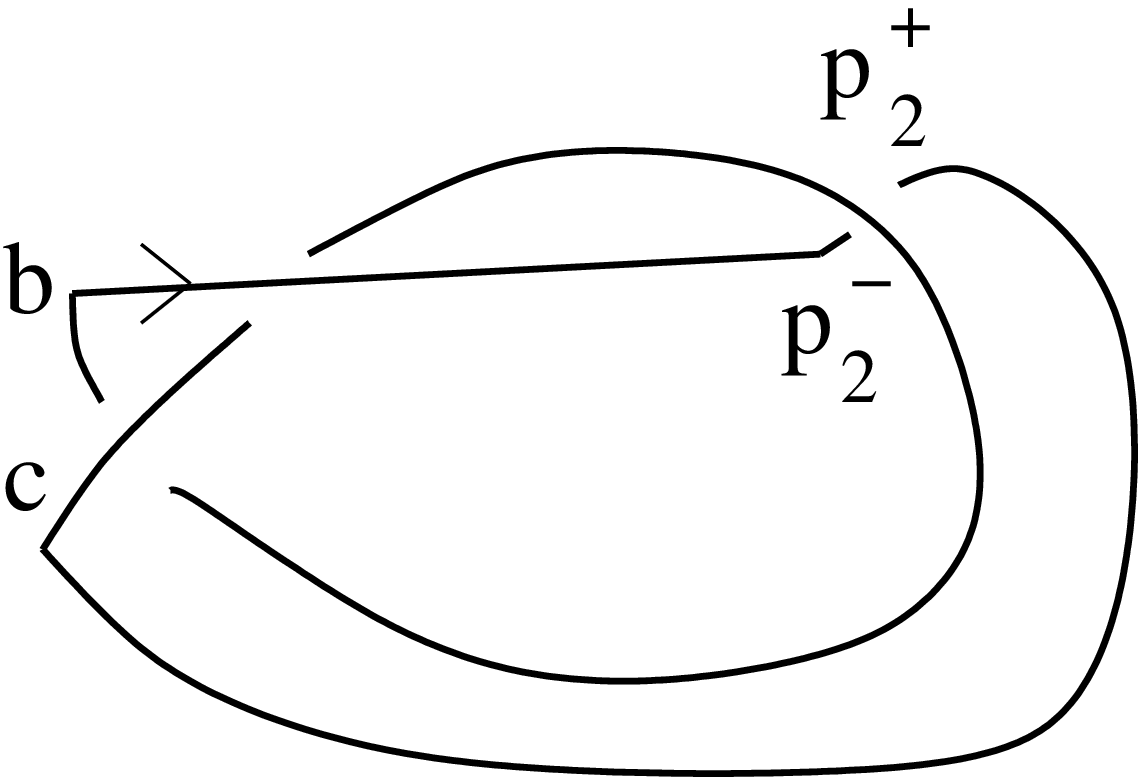,height=3.4cm}}
\centerline{Fig. 3.}
\ \\
This completes the proof of Theorem 1.
in the case of a knot. The case of a link is analogous. The only
difference is that we have to consider an innermost 1-gon or 0-gon
and reduce diagram to the right handed trefoil or Hopf link (with,
possibly, some additional trivial components). 
\end{proof}

Theorem 1 is stronger than that of Rudolph or Murasugi because
there are positive knots which are neither alternating nor have a
presentation as positive braids.

\begin{corollary}
  A nontrivial positive link is neither slice link nor
  amphicheiral link.
\end{corollary}
\begin{proof}
  Amphicheiral links and slice links have signature equal to 0
  \cite{M-1}.
\end{proof}
\begin{corollary}[Murasugi]\label{Cor3}
The following inequalities hold for the Jones polynomial of a
non-split link:
\begin{itemize}
    \item [(a)] $d_{min}V_L(t) >0$ for a nontrivial positive link
    $L$,
    \item [(b)]$d_{max}V_L(t) <0$ for a nontrivial negative link
    $L$,
\end{itemize}
where $d_{max}$ (resp. $d_{min}$) denotes the highest (resp.
lowest) power of $t$ in $V_L(t)$.
\end{corollary}
\begin{proof}
  By Murasugi \cite{M-5}, Theorem 13.3, the following holds for any
  non-split link diagram $\widetilde{L}$ of a link L: $\left\{%
\begin{array}{ll}
    d_{max}V_L(t) \leq c_+((\tilde L))- \frac{1}{2} \sigma(L), & \hbox{} \\
     d_{min}V_L(t) \geq -c_-((\tilde L))- \frac{1}{2} \sigma(L), & \hbox{} \\
\end{array}%
\right.    \ \\  
$ where $c_+$ (resp. $c_-$) is the number of positive
(resp. negative) crossings of $\tilde L$. Now Corollary
\ref{Cor3} follows from Theorem \ref{thm1}.
\end{proof}
Corollary \ref{Cor3} was first proven by Murasugi [\cite{M-3},
Theorem 2.1] in implicit form. A different proof has been found by
Traczyk. 

Theorem \ref{thm1} can be extended to other Tristram--Levine
signatures as long as it holds for the Hopf link and the trefoil
knot.\\
We use the notation of \cite{G}(see also \cite{P}). We assume also
(without loss of generality) that $|1-\xi|=1$ in the
Tristram--Levine signature $\sigma_{\xi}$.

\begin{theorem}
If $L$ is a nontrivial positive link then for $Re \xi < \frac{1}{2}$, \ \ 
$\sigma_{\xi} <0$.
\end{theorem}
\begin{proof}
  For $Re \xi <1/2$, $\sigma_{\xi}$ is negative for the right
  handed trefoil knot and Hopf link. Furthermore, by \cite{P-T}
(see also [\cite{P}, Lemma 4.13(b)] for $Re \xi <1$,
$\sigma_{\xi}(L_+) \leq \sigma_{\xi}(L_-)$ so the proof of Theorem
\ref{thm1} can be repeated without changes here too.
\end{proof}
\begin{conjecture}
  If a nontrivial link has a diagram with at most one negative
  crossing then the link has negative signature.\footnote{Added for e-print:\ This conjecture 
with its generalizations was proved in a joint paper with K.~Taniyama \cite{P-T}.}
\end{conjecture}

I have been informed, after completing this manuscript that the
Theorem \ref{thm1} has been proven independently by P. Traczyk
("Non-trivial negative links have positive signature", preprint,
Summer 1987) and, in the case of knots by R. Gompf and T. Cochran
(``Applications of Donaldson's theorems to classical knot
concordance. Homology 3-spheres and property P", preprint 1987).

\centerline{Department of Mathematics}
\centerline{University of Toronto}
\centerline{Toronto, Canada}
\centerline{M5S 1A1}
\centerline{(and Warsaw University)}
\ \\
This e-print is based on the preprint written in  Toronto in September 1987; most likely 
it is the version I submitted to Bulletin Polish Acad. Sci.: Math.


\begin{thebibliography}{99}
\bibitem[B-W]{B-W} J. Birman, R.F. Williams, Knotted orbits in
dynamical systems- I: Lorentz's equations, {\it Topology}, 22(1), 1983,
47-82.

\bibitem[G]{G} C. McA. Gordon, some aspects of classical knot
theory, In: Knot theory, Lect. Notes in Math. 685, 1978, 1-60.

\bibitem [M-1]{M-1} K.~Murasugi, On a certain numerical invariant
of link types, {\it Trans. Amer. Math. Soc.} 117, 1965, 387-422.

 \bibitem [M-2]{M-2} K.~Murasugi, On the signature of links,
 {\it Topology} 9, 1981, 283-298.

  \bibitem [M-3]{M-3} K.~Murasugi, Jones polynomial of alternating
  links, {\it Trans. Amer. Math. Soc.} 295(1), 1986, 147-174.

\bibitem [M-4]{M-4}
K.~Murasugi, Jones polynomial and classical conjectures in knot
theory, II, {\it Math. Proc. Camb. Phil. Soc.}, to appear. 
(Added for e-print:\ 102, 1987, 317-318.)

 \bibitem [M-5]{M-5} K.~Murasugi, On invariants of graphs with
 applications to knot theory, preprint 1987.\\
(Added for e-print:\ {\it Trans. Amer. Math. Soc.}, 314, 1989, 1-49.)

\bibitem [P]{P}
J.~H.~Przytycki, Survey on recent invariants  in classical
knot theory,  preprint, Warsaw University, 1986. 
(Added for e-print:\ Warsaw University, Preprints 6,8,9; Warszawa, 1986;\
e-print:\  {\tt http://front.math.ucdavis.edu/0810.4191})

\bibitem [P-T]{P-T}
J.~H.~Przytycki, P.~Traczyk, Conway algebras and skein equivalence
of links, preprint Warsaw 1985. 
(Added for e-print:\ part of the preprint was published in 
{\it Proc. Amer. Math. Soc.},  100(4), 1987, 744-748.)

\bibitem[R]{R} L. Rudolph, Non-trivial positive braids have
positive signature, {\it Topology} 21, 1982, 325-327.
\\ \ \\
{\bf Added for e-print:}\ 
\bibitem [C-G]{C-G} 
T.~Cochran, E.~Gompf, Applications of Donaldson's theorems to
classical knot concordance, homology $3$-spheres and property $P$,
{\it Topology} 27(4), 1988, 495--512.

\bibitem [P-T]{P-T} 
J.~H.~Przytycki, K.~Taniyama, 
 Almost positive links have negative signature, preprint 1991;\
e-print: \  {\tt arXiv:0904.4130 }

\bibitem
[T-1]{T-1}
K.~Taniyama, A partial order of knots, {\it Tokyo J. Math.}
      12(1), 1989, 205-229.

\bibitem
[T-2]{T-2}
K.~Taniyama, A partial order of links, {\it Tokyo J. Math.}
      12(2), 1989, 475-484.

\bibitem [T]{T}
P.~Traczyk,  Nontrivial negative links have positive signature.
      {\it Manuscripta Math.} 61(3), 1988, 279--284.

\end{thebibliography}
\end{document}